\documentclass[leqno]{amsart}
\usepackage{amscd,amssymb,amsmath}

\usepackage{tikz}
\usepackage{color}
\usepackage{float}
\usepackage{graphicx}
\usepackage{vntex}
\usepackage[brazil]{babel}    

\usepackage{amsmath,amsfonts,amsthm,amssymb}
\usepackage[all]{xy}     
\numberwithin{equation}{section}

\theoremstyle{definition}
\newtheorem{theorem}[equation]{Teorema}

\newtheorem{definition}[equation]{Defini\c c\~ao}

\numberwithin{equation}{section}
\theoremstyle{definition}

\newcommand{\R}{\mathbb{R}}
\newcommand{\Sp}{\mathbb{S}}
\newcommand{\T}{\mathbb{T}}
\newcommand{\PP}{\mathbb{P}}

\begin{document}

\title[A caracter\'\i stica de Descartes-Poincar\'e] {A caracter\'istica de Euler-Poincar\'e}

\author{Jean-Paul \textsc{Brasselet} e Nguy\~{\^e}n Th\d{i} B\'ich Th\h{u}y}
\address{CNRS I2M e Aix-Marseille Universit\'e, Marseille, France.}
\email{jean-paul.brasselet@univ-amu.fr}

\address{UNESP, Universidade Estadual Paulista, ``J\'ulio de Mesquita Filho'', S\~ao Jos\'e do Rio Preto, Brasil}
\email{bichthuy@ibilce.unesp.br}

\bigskip



\begin{abstract}
D\^e a C\'esar o que \'e de C\'esar e a Descartes o que \'e de Descartes.

\end{abstract}

\maketitle


\section{Introdu\c c\~ao}
Netse artigo, vamos introduzir a caracter\'\i stica de Euler-Poincar\'e, de uma  maneira elementar, assim
que historicamente. Veremos tamb\'em porque deveriar chamar-se de caracter\'\i stica de Descartes-Poincar\'e.

Todos espa\c cos ser\~ao compactos, sem bordo. Vamos trabalhar essencialmente com superf\'icies suaves e orient\'aveis. Por\'em, olharemos tamb\'em superf\'icies singulares ou n\~ao orient\'aveis.  Daremos tamb\'em alguns resultados que estendem-se no caso de dimens\~ao maior do que 2.

O artigo de Elon Lima \cite{Li2} sobre a caracter\'\i stica de Euler-Poincar\'e, publicado no primeriro 
n\'umero da Revista Matem\'atica Universit\'aria, fornece defini\c c\~oes equivalentes desta invariante. 

\section{O Teorema de Descartes-Euler}

O nome de ``caracter\'\i stica de Euler-Poincar\'e" vem do an\'uncio de Leonhard Euler, numa carta a um amigo, Goldbach (1751) do resultado seguinte:

\begin{theorem}
Seja $K$ uma triangula\c c\~ao (Defini\c c\~ao \ref{triangulation}) 
da esfera $\Sp^2$, com  n\'umeros $n_0$ de v\'ertices, $n_1$ de segmentos e $n_2$ de tri\^angulos, ent\~ao temos
\begin{equation}\label{Euler}
n_0 - n_1 + n_2 = 2.
\end{equation}
\end{theorem}

\begin{figure}[h]
\hskip 20 pt
\scalebox{1.20}{\includegraphics{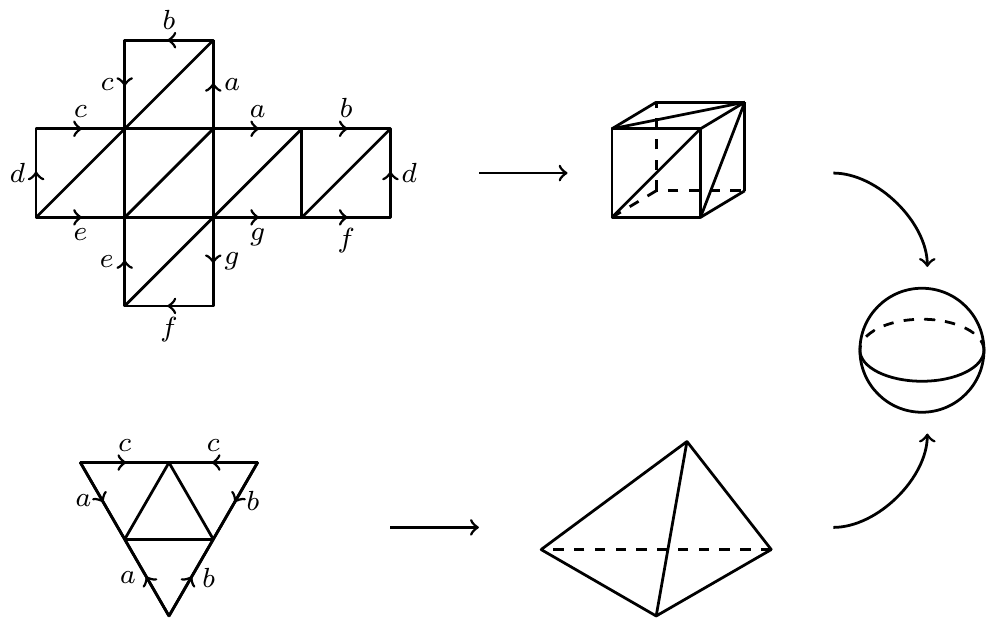}}
\caption{Triangula\c c\~oes da esfera.}\label{triangsphere}
\end{figure}

Na figura \ref{triangsphere}, temos dois exemplos de triangula\c c\~oes da esfera $\Sp^2$ nos quais se pode verificar a 
rela\c c\~ao (\ref{Euler}). Pela triangula\c c\~ao do cubo, temos 
$$n_0 - n_1 + n_2 = 8 - 18 + 12 =2.$$
Pelo tetraedro, temos 
$$n_0 - n_1 + n_2 = 4 -6 + 4 =2.$$

De fato, Euler nunca deu uma prova correta deste resultado. Al\'em disso, 
o resultado j\'a era provado por Descartes (1639). 
Descartes morreu em 1650. 
O editor de Descartes comunicou os manuscritos de Descartes com Leibnitz, que ficou em Paris de 1672 at\'e 1676.  Leibnitz tinha vontade de publicar o manuscrito "De solidorum elementis" de Descartes, mas nunca concretizou este projeto. 
Em 1703, E. de Jonqui\`eres publicou uma ``Notes aux Comptes Rendus de l'Acad\'emie des Sciences"
\cite{dJ1} onde ele contou a hist\'oria do manuscrito (seguido Foucher de Carell \cite{Fo}) e mostrou que Descartes j\'a conheceu a rela\c c\~ao (\ref{Euler}) chamada de 
Euler. 
Alguns autores dizem mesmo que Archimedes j\'a conheceu o resultado. 

Houve muitas provas da f\'ormula. 
O site web \cite{Epp} fornece 19 provas poss\'iveis, incluindo as provas de Descartes (usando \^angulos) e de Cauchy. 
De fato, depois de Descartes, uma das primeiras provas e mais bonita foi dada por Cauchy (1811). 
Vamos dar aqui uma id\'eia desta prova. 
Antes disso, precisamos providenciar algumas defini\c c\~oes. 

Vamos considerar poliedros particulares, isto \'e, poliedros munidos de uma  triangula\c c\~ao no sentido seguinte: 
Um poliedro ser\'a a uni\~ao (num espa\c co euclidiano $\R^n$) de um conjunto finito de 
``simplexos'' : tri\^angulos (simplexos de dimens\~ao 2, denotados por $\sigma^2$), segmentos (simplexos de dimens\~ao 1, denotados por $\sigma^1$) e v\'ertices (simplexos de dimens\~ao 0, denotados por $\sigma^0$), tais que cada 
face de um simplexo \'e um elemento 
do poliedro e para todo par $\sigma_i, \sigma_j$ de simplexos, 
a intersec\c c\~ao \'e seja vazia, seja uma face comum de $\sigma_i$ e $ \sigma_j$. A uni\~ ao dos 
simplexos do poliedro \'e um subespa\c co topol\'ogico compacto de $\R^n$ denotado por $\vert K \vert$.

Os poliedros ser\~ao vizualizados seja num espa\c co euclidiano (em geral $\R^3$ - veja os exemplos do cubo ou do tetraedro na figura \ref{triangsphere}),
 ou seja sub a forma de representa\c c\~ao plana. Neste caso, o bordo da figura consiste em segmentos 
que s\~ao identificados com uma orienta\c c\~ao dada (veja as representa\c c\~oes planas 
do cubo e do tetraedro na figura \ref{triangsphere}). 

\begin{figure}[h!] 
\label{AA}
\begin{tikzpicture} 
  \draw (1,0) node[below]{$B$}  -- (0,1) node[left]{$A$}  -- (1,2) node[above]{$A$} -- (1,0);  
\draw (4,0) node[below]{$B$}  -- (3,1.5) node[left]{$A$}  -- (4,2) node[above]{$C$} -- (4,0); 
\draw (4,0)  -- (4,2)  -- (6, 3) node[above]{$A$} -- (4,0); 
\end{tikzpicture}
\caption{Tri\^angulos proibidos.}
\end{figure}
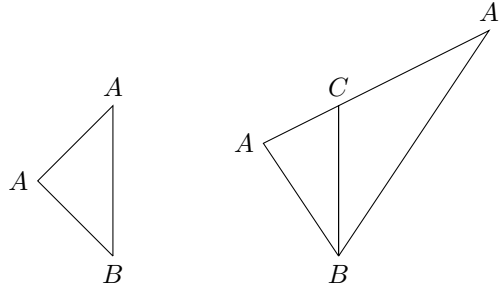

Uma observa\c c\~ao importante \'e que numa representa\c c\~ao plana dum poliedro de dimens\~ao 2, 
s\~ao proibidos os tri\^angulos seguintes (figura \ref{AA}): O primeiro daria um segmento $AB$ no poliedro e o segundo um tri\^angulo $ABC$ somente, 
e n\~ao dois tri\^angulos como esperado.

\medskip

\medskip

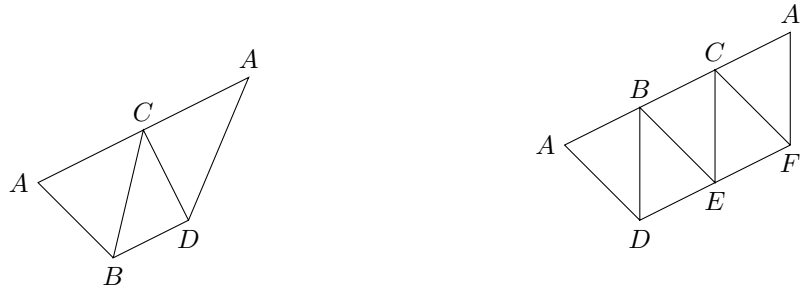
\begin{figure} [h!]
 \label{AB}
\begin{tikzpicture} 
  \draw (1, 0) node[below]{$B$}  -- (0, 1) node[left]{$A$}  -- (1.4,1.7) node[above]{$C$} -- 
(2.8,2.4) node[above]{$A$} -- (2,0.5) node[below]{$D$} -- (1, 0) --(1.4,1.7) -- (2, 0.5);  
 
\draw (8, 0.5) node[below]{$D$}  -- (7,1.5) node[left]{$A$}  -- (8,2) node[above]{$B$} -- 
(9,2.5) node[above]{$C$} -- (10,3) node[above]{$A$} -- (10, 1.5) node[below]{$F$}-- 
(9,1) node[below]{$E$}
 -- (8, 0.5) -- (8,2)  --  (9,1) -- (9, 2.5) -- (10, 1.5);  
\end{tikzpicture}

\caption{Tri\^angulos permitidos.}
\end{figure}

 Os desenhos na Figura \ref{AB} s\~ao permitidos.

\medskip



\begin{definition}\label{triangulation}
 Uma triangula\c c\~ao da esfera \'e o dado de um poliedro (triangulado) $K$ e de um
homeomorfismo $h:\vert K \vert \to X$.
\end{definition} 

Agora, vamos \`a prova do Cauchy. A primeira obseva\c c\~ao \'e que, se tiramos um ponto
$\{ pt \}$ da esfera,  ent\~ao o resto $\Sp^2 \setminus \{ pt \}$ \'e contract\'ivel num ponto. 

Vamos fazer a conta $n_0 - n_1 + n_2 $ pelo poliedro $K$ que representa uma triangula\c c\~ao da esfera $\Sp^2$. Numa primeira etapa, tiramos um tri\^angulo (aberto) $\sigma^2$, ent\~ao j\'a temos $1$ na conta de $n_2$. Vamos chamar de $B$ o buraco assim formado (veja figura \ref{figura01}). 

A segunda etapa consiste em duas opera\c c\~oes que vamos descrever:

O bordo de $B$ est\'a constituido de
tr\^es segmentos $\sigma^1_i$ com a seguinte propriedade : 
Cada $\sigma^1_i$ \'e uma face de um tri\^angulo $\sigma^2_i$ que tem somente o segmento $\sigma^1_i$ 
como face comum com o buraco $B$.

\begin{figure} [h!]
\begin{tikzpicture} 
  \draw (1,0)   --  (0,2) -- (1,4) -- (2,5) -- (4, 4) -- (5, 3) -- (4,1) -- (2.5, 0.25) --(2,1)  -- (1.5,2) -- (2,3) -- (3,2) -- (4,2) ;  

\draw (1, 0) -- (2.5, 0.25); 
\draw (1, 0) -- (2, 1);
\draw (1, 0) -- (1.5, 2);

\draw (0,2) -- (1.5, 2);
\draw (1,4) -- (1.5,2);
\draw (1,4)  -- (2,3);

\draw (2,5) -- (2,3);

\draw (4,4) -- (2,3);
\draw (4,4) -- (3,2);
\draw (4,4) -- (4,2);

\draw (5,3) -- (4,2);

\draw (4,1) -- (4,2);
\draw (4,1) -- (3,2);
\draw (4,1) -- (2,1);

\draw (2.5,0.25) -- (1,0);

\draw (2,1) -- (3,2);

\draw (1.5,2) --(3,2);

\node (a) at (2.2, 1.6) {$B$};
\node (b) at (2.2, 2.2) {$\sigma_i^1$}; 
\node (b) at (2.1, 2.6) {$\sigma_i^2$};

\filldraw[fill=blue!20] (2,1) -- (1.5, 2) -- (3,2);
\node (a) at (2.2, 1.6) {$B$};
 
\end{tikzpicture}
\caption{A opera\c c\~ao ``I''.}\label{figura01}
\end{figure}


A opera\c c\~ao ``I'' consiste a tirar tal tri\^angulo $\sigma^2_i$ junto com o segmento $\sigma^1_i$ 
correspondente e assim o buraco se estende. 

Qual \'e a consequ\^encia da opera\c c\~ao ``I'' sobre a soma $n_0 - n_1 + n_2$ ?

Claramente, a opera\c c\~ao ``I'' n\~ao muda a soma $n_0 - n_1 + n_2$ porque ela diminui 
$n_1$ e $n_2$ de $1$, ent\~ao $n_0 - n_1 + n_2$ fica igual.

Continuando o processo, chegaremos \`a situa\c c\~ao da figura \ref{figura03}, onde um tri\^angulo 
(aqui $\tau^2$) tem duas faces (segmentos) comuns $\tau^1_i$ e $\tau^1_j$  com o buraco.

\begin{figure} [h!]
\begin{tikzpicture} 
  \draw (1,0)   --  (0,2) -- (1,4) -- (2,5) -- (4, 4) -- (5, 3) -- (4,1)  --(2,1)  -- (1.5,2) -- (2,3) -- (3,2) -- (4,2) ;  

\draw (1, 0) -- (2, 1);
\draw (1, 0) -- (1.5, 2);

\draw (0,2) -- (1.5, 2);
\draw (1,4) -- (1.5,2);
\draw (1,4)  -- (2,3);

\draw (2,5) -- (2,3);

\draw (4,4) -- (2,3);
\draw (4,4) -- (3,2);
\draw (4,4) -- (4,2);

\draw (5,3) -- (4,2);

\draw (4,1) -- (4,2);
\draw (4,1) -- (3,2);
\draw (4,1) -- (2,1);


\draw (2,1) -- (3,2);

\draw (1.5,2) --(3,2);

\draw (0,0) -- (2, -1);
\draw (2,-1) --(1,0);
\draw (2,-1) -- (4,1);
\draw (0,0) -- (1,0);
\draw (0,0) -- (0,2);
\draw (1,0) -- (4,1);

\draw[red, line width=2pt] (2,1) -- (1.5,2);
\draw[red, line width=2pt] (2,1) -- (3,2);

\filldraw[fill=red!20] (2,1) -- (1.5, 2);

\filldraw[fill=blue!20] (1,0) -- (2, 1) -- (1.5,2);
\filldraw[fill=blue!20] (1,0) -- (2, 1) -- (4,1);
\filldraw[fill=blue!20] (4,1) -- (2, 1) -- (3,2);

\node (a) at (2.2, 1.7) {$\tau^2$};
\node (b) at (1.6, 1.4) {$\tau_i^1$}; 
\node (c) at (2.8, 1.5) {$\tau_j^2$};
\node (d) at (2.2, 0.8) {$\tau^0$};


\end{tikzpicture}
\caption{A opera\c c\~ao ``II''. }\label{figura03}
\end{figure}


A opera\c c\~ao ``II'' consiste a tirar tal tri\^angulo $\tau^2$ junto com os dois segmentos 
$\tau^1_i$ e $\tau^1_j$ e com o v\'ertice $\tau^0$ que \'e face comum de $\tau^1_i$ e $\tau^1_j$. 
Assim, o buraco se estende. 

Qual \'e a consequ\^encia da opera\c c\~ao ``II'' sobre a soma $n_0 - n_1 + n_2$ ?

Claramente, a opera\c c\~ao ``II'' n\~ao muda a soma $n_0 - n_1 + n_2$ porque ela diminui 
$n_0$ e $n_2$ de $1$ e diminui $n_1$ de 2, ent\~ao $n_0 - n_1 + n_2$ fica igual.

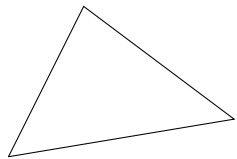
\begin{figure} [h!]
\begin{tikzpicture} 
  \draw (0,0)   --  (1,2) -- (3,0.5) -- (0,0)  ;  

\end{tikzpicture}
\caption{O \'ultimo tri\^angulo.}\label{figura04}
\end{figure}

Continuando este processo, pelas opera\c c\~oes ``I'' e ``II''  temos de cuidar que o bordo do buraco seja sempre uma curva ``simples'' (homeomorfa a um c\'irculo) ao fim de ficar com uma figura conexa. 
Sabendo que $\Sp^2 \setminus B$ \'e 
contract\'ivel, chegamos a tirar todos tri\^angulos ao menos do \'ultimo que fica sozinho 
(veja figura \ref{figura04}).

Nesta \'ultima etapa, a conta de $n_0 - n_1 + n_2$ pelo \'ultimo tri\^angulo d\'a :
$$n_0 - n_1 + n_2 = 3 - 3 + 1 = 1.$$
Mas, n\~ao esquecemos que, na primeira etapa, tiramos um tri\^angulo, que contribui a $n_2$. 
Ent\~ao o total est\'a bem 
$$n_0 - n_1 + n_2 = 2,$$
o que prova o resultado.

Observamos que o resultado vale para poliedros n\~ao necessariamente triangulados. Mais precisamente, temos o teorema :

\begin{theorem}
Seja $P$ um poliedro homeomorfa \`a esfera $\Sp^2$. Sejam $n_0$ o n\'umero de v\'ertices, $n_1$ o n\'umero de arestas (segmentos) e $n_2$ o n\'umero de pol\'igonos de dimens\~ao 2, ent\~ao
$$n_0 - n_1 + n_2 = 2.$$
\end{theorem}

A prova deste resultado consiste a ``triangular" os pol\'igonos e mostrar que a soma $n_0 - n_1 + n_2$ 
n\~ao muda pelo processo. Suponhamos que os v\'ertices dum pol\'igono s\~ao os pontos 
$a_1, a_2, a_3, a_4, \ldots , a_k$. 

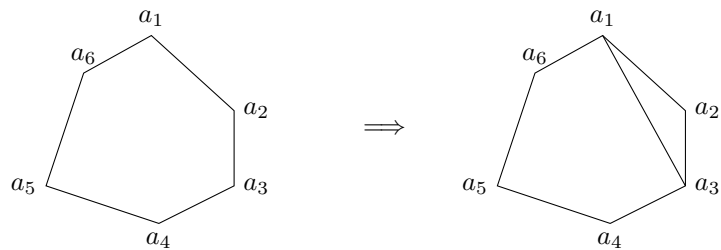
\begin{figure}  [h!]
\begin{tikzpicture} 
  \draw (1.5,0) node[below]{$a_4$}   --  (0,0.5) node[left]{$a_5$} -- (0.5, 2) node[above]{$a_6$} -- ( 1.4, 2.5) node[above]{$a_1$} -- (2.5, 1.5) node[right]{$a_2$} -- (2.5,0.5) node[right]{$a_3$}  -- (1.5,0) ;  

\node (a) at (4.5, 1.25) {$\Longrightarrow$};

  \draw (7.5,0) node[below]{$a_4$}   --  (6,0.5) node[left]{$a_5$} -- (6.5, 2) node[above]{$a_6$} -- ( 7.4, 2.5) node[above]{$a_1$} -- (8.5, 1.5) node[right]{$a_2$} -- (8.5,0.5) node[right]{$a_3$} -- (7.5,0);

\draw ( 7.4, 2.5) -- (8.5, 0.5);

\end{tikzpicture}
\caption{Triangula\c c\~ao dum pol\'igono.}\label{AC}
\end{figure}


Se juntamos os pontos $a_1$ e $a_3$ por um segmento, para formar um tri\^angulo $(a_1, a_2, a_3)$, temos uma subdivis\~ao do pol\'igono $(a_1, a_2, a_3, a_4, \ldots , a_k)$ em um 
tri\^angulo $(a_1, a_2, a_3)$ e um novo pol\'igono $(a_1, a_3, a_4, \ldots , a_k)$. O n\'umero de 
pol\'igonos sube de $1$ e tamb\'em o n\'umero de arestas (temos uma nova aresta $(a_1,a_3)$). Ent\~ao a soma $n_0 - n_1 + n_2$ n\~ao muda. Continuando este processo duma maneira \'obvia, obtemos um poliedro triangulado e terminamos com a prova anterior de Cauchy. 

\begin{figure}[h]
\hskip 20 pt
\scalebox{0.30}{\includegraphics{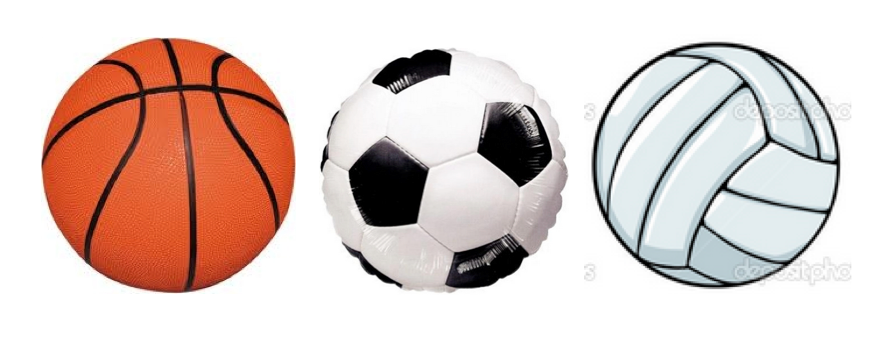}}
\caption{Poliedros de bal\~oes de esportes.}\label{esportes}
\end{figure}

Assim, o desenho dum bal\~ao de futebol 
(figura \ref{esportes}) \'e  um exemplo de poliedro (n\~ao triangulado) que satisfaz 
a f\'ormula. Isso vale para os poliedros desenhados nos outros bal\~oes de esporte !

\section{O toro e demais superf\'icies}

\subsection{O Toro $\T$}

Podemos calcular a caracter\'\i stica de Euler-Poincar\'e do toro da mesma maneira 
do que a a caracter\'\i stica de Euler-Poincar\'e da esfera. 

\begin{figure}[h]
\hskip 20 pt
\scalebox{1.40}{\includegraphics{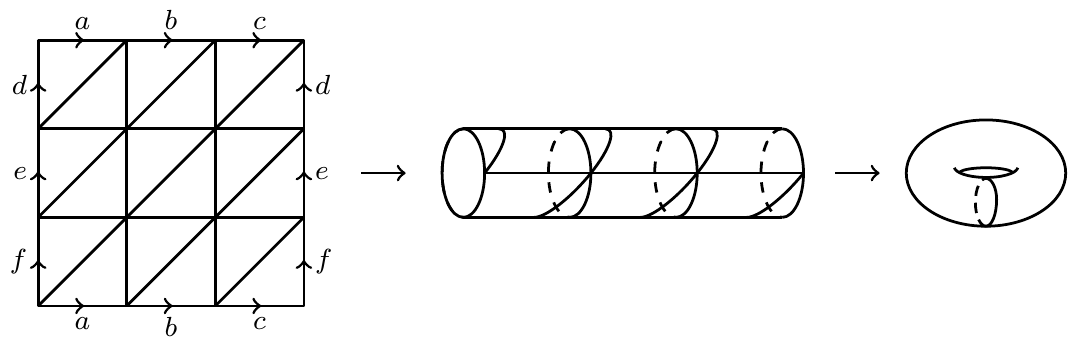}}
\caption{Representa{\c c}\~ao plana do toro.}\label{toro02}
\end{figure}

A representa{\c c}\~ao plana do toro  \'e desenhada na figura \ref{toro02}.

Assim temos $$n_0 - n_1 + n_2 =  9 - 27 + 18= 0.$$

\begin{figure}[h]
\hskip 20 pt
\scalebox{0.40}{\includegraphics{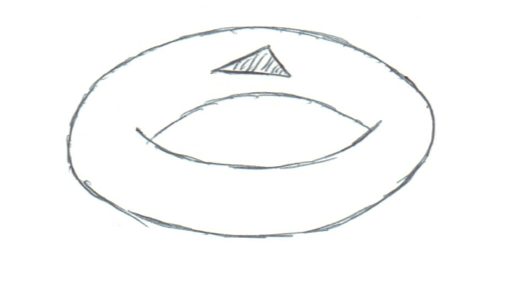}}
\caption{Tiramos um tri\^angulo do toro.}\label{figura05}
\end{figure}

Para mostrar que esta quantidade n\~ao depende da triangula\c c\~ao do toro, 
podemos usar a mesma t\'ecnica do que pela esfera, seguindo a prova de Cauchy :

Numa primeira etapa, tiramos um tri\^angulo (aberto) $\sigma^2$, ent\~ao j\'a temos $1$ na conta de $n_2$ (Figura \ref{figura05}).

%

\begin{figure}[h]
\hskip 20 pt
\scalebox{0.40}{\includegraphics{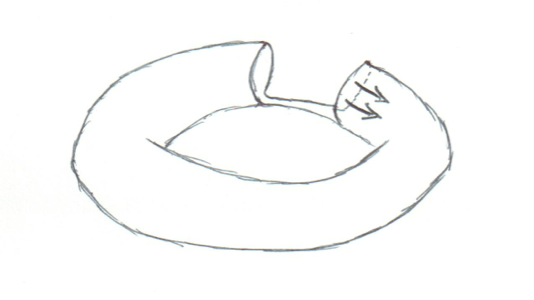}}
\caption{Extens\~ao do buraco - I.}\label{figura06}
\end{figure}

\begin{figure}[h]
\hskip 20 pt
\scalebox{0.40}{\includegraphics{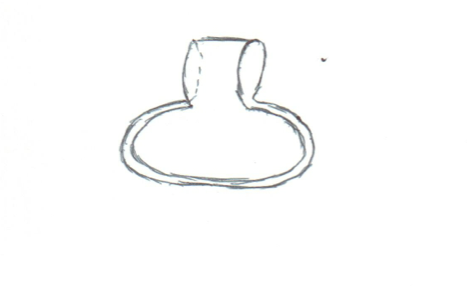}}
\caption{Extens\~ao do buraco - II.}\label{figura07}
\end{figure}

O buraco assim efetuado se estende, usando as opera\c c\~oes ``I'' e ``II'' (Figuras \ref{figura06} 
e \ref{figura07}) efetuando a extens\~ao no sentido das flechas. 

\begin{figure}[h]
\hskip 20 pt
\scalebox{0.40}{\includegraphics{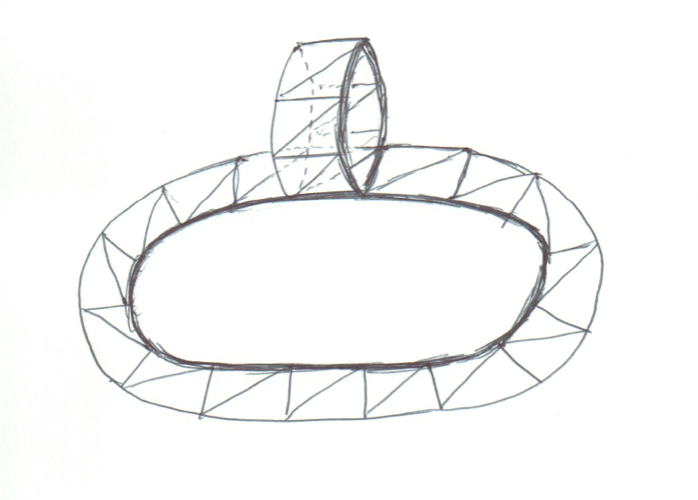}}
\caption{As faixas.}\label{figura08}
\end{figure}

Agora, temos a figura de duas faixas, coladas 
e trianguladas como indicado na figura \ref{figura08}.

\begin{figure}[h]
\hskip 20 pt
\scalebox{0.8}{\includegraphics{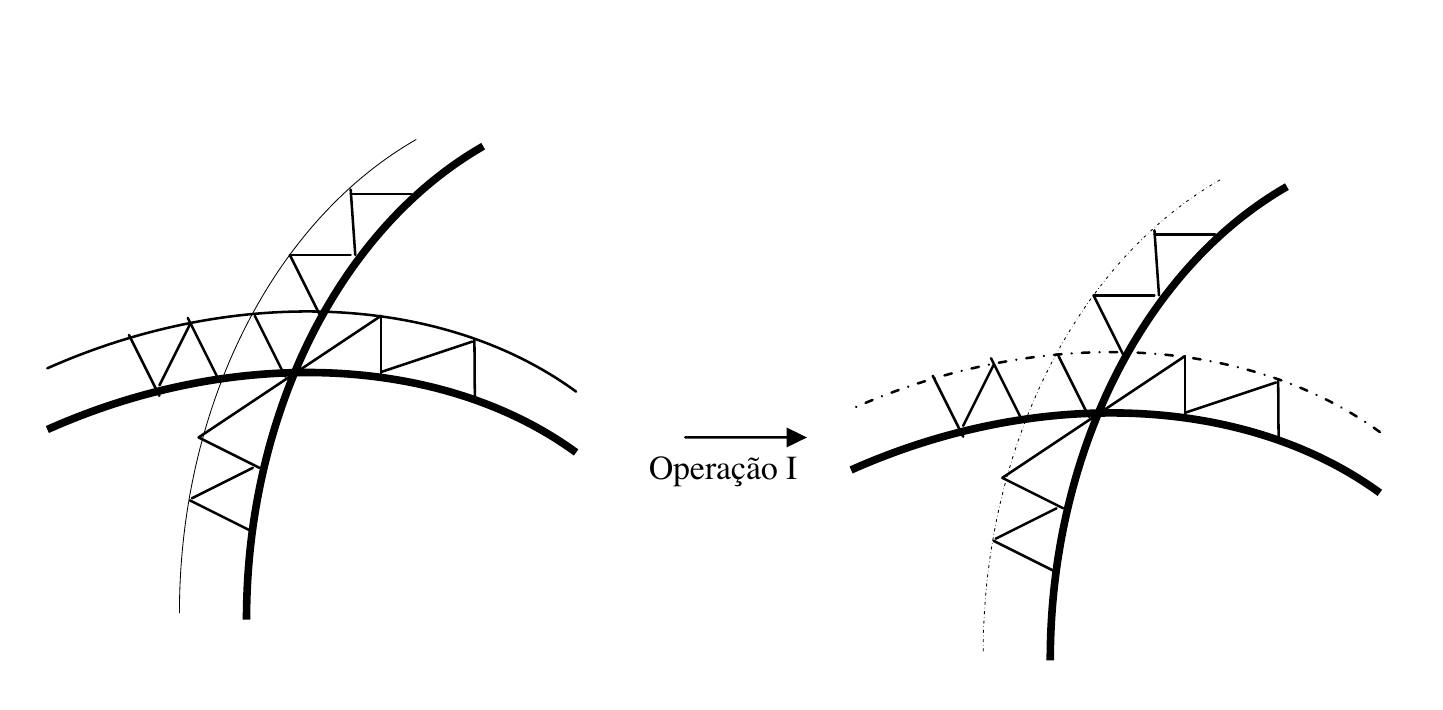}}
\caption{Opera{\c c}\~ao ``I'' perto do cruzamento das faixas.}\label{figura09}
\end{figure}

\begin{figure}[h]
\hskip 20 pt
\scalebox{0.80}{\includegraphics{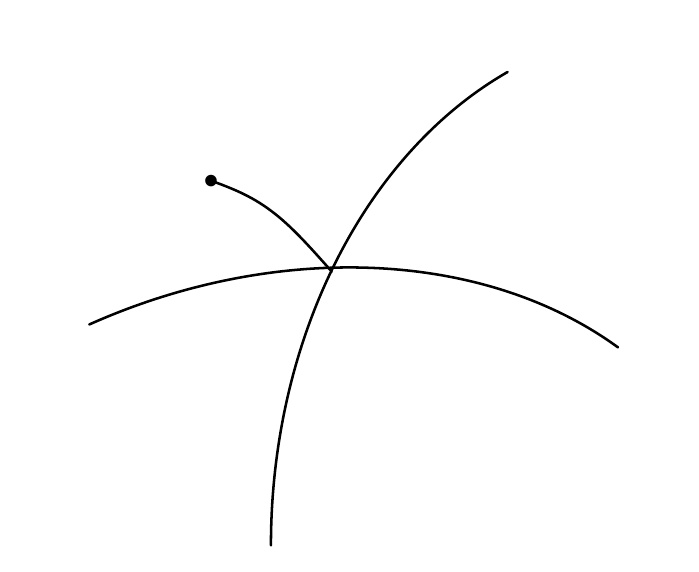}}
\caption{A figura depois das opera\c c\~oes ``I'' e ``II''.}\label{figura10}
\end{figure}

Fora do cruzamento das faixas, o uso das opera\c c\~oes ``I'' e ``II''  \'e \'obvio.
No cruzamento, se usa a opera\c c\~ao ``I'' (Figura \ref{figura09} e depois 
a opera\c c\~ao ``II''). Obtemos o desenho da Figura \ref{figura10}.

Tirando o segmento que fica e o v\'ertice correspondente n\~ao muda a soma 
$n_0 - n_1 + n_2$. 

Agora, temos dois c\'irculos com um ponto comum
(Figura \ref{figura11}). Qualquer triangula\c c\~ao desta figura
fornece um n\'umero $k$ de v\'ertices e $(k+1)$ de segmentos. Ent\~ao, para esta figura,
$n_0 - n_1 + n_2 = -1$. Mas, n\~ao esquecemos que, na primeira etapa, tiramos um tri\^angulo, que contribui a $n_2$. 
Ent\~ao o total est\'a bem 
$$n_0 - n_1 + n_2 = 0,$$
qualquer que seja a triangula\c c\~ao do toro. 

\begin{figure}[h]
\begin{tikzpicture} 
\draw (2,2) ellipse (2cm and 1cm);
\draw (2,4) ellipse (0.5cm and 1cm);
\node at (2,3) {$\bullet$};
\end{tikzpicture} 
\caption{A figura final.}\label{figura11}
\end{figure}


\vfill\break
Na seguinte, denotaremos por $\chi({\mathcal S})$ a caracter\'\i stica de Descartes-Poincar\'e de uma superf\'icie $\mathcal S$. 
\medskip

Usaremos tamb\'em no\c c\~oes de g\^enero e de soma direta. 

O {\it g\^enero} $g$ de uma superf\'icie (orient\'avel ou n\~ao) 
\'e o n\'umero m\'aximo de c\'irculos que se pode desenhar sobre  a superf\'icie sem a desconectar.

A soma conexa de duas superficies significa que tiramos um disco de cada uma das superf\'icies e que colamos as supreficies ao longo dos c\'irculos obtidos (bordos dos buracos que fizemos) (veja figura \ref{AE}).

\begin{figure}[h]
\hskip 20 pt
\scalebox{0.70}{\includegraphics{AE.jpg}}
\caption{Soma conexa.}\label{AE}
\end{figure}

\section{Superf\'icies orient\'aveis e n\~ao orient\'aveis}

\subsection{Orienta{\c c}\~ao}
Denotamos por $a_i$ os v\'ertices dos simplexos (v\'ertices, segmentos e tri\^angulos). 
Quando definimos uma ordem pelos v\'ertices de segmentos (tais que $(a_0, a_1)$) 
e de tri\^angulos (tais que $(a_0, a_1,a_2)$), 
denotamos os respectivamente por $[a_0, a_1]$ ou $[a_1, a_0]$ no caso do segmento e, por exemplo,  
$[a_0, a_1, a_2]$ ou  $[a_1, a_2, a_0]$, etc... no caso do tri\^angulo.

Dar uma ordem $[a_0, a_1, a_2]$ dos v\'ertices de um tri\^angulo, defina uma orienta{\c c}\~ao do tri\^angulo (Figura \ref{figura12a}) :

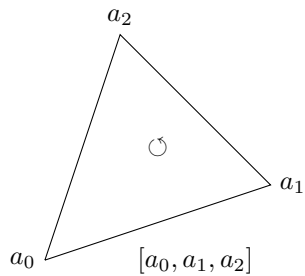
\begin{figure}[h]
\begin{tikzpicture} 
  \draw (0,0) node[left]{$a_0$}   --  (1,3) node[above]{$a_2$} -- (3, 1) node[right]{$a_1$} -- (0,0);
\node at (1.5, 1.5) {$\circlearrowleft$}; 
\node at (2,0) {$[a_0, a_1, a_2 ]$};
\end{tikzpicture}
\caption{Orienta{\c c}\~ao do tri\^angulo.}\label{figura12a}
\end{figure}



Fazendo uma permuta{\c c}\~ao de ordens (por exemplo $[a_0, a_2, a_1]$), defina a orienta{\c c}\~ao oposta.

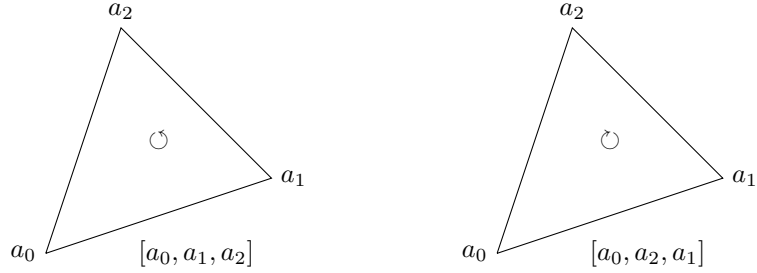
\begin{figure}[h]
\begin{tikzpicture} 
  \draw (0,0) node[left]{$a_0$}   --  (1,3) node[above]{$a_2$} -- (3, 1) node[right]{$a_1$} -- (0,0);
\node at (1.5, 1.5) {$\circlearrowleft$}; 
\node at (2,0) {$[a_0, a_1, a_2 ]$};

\draw (6,0) node[left]{$a_0$}   --  (7,3) node[above]{$a_2$} -- (9, 1) node[right]{$a_1$} -- (6,0);
\node at (7.5, 1.5) {$\circlearrowright$}; 
\node at (8,0) {$[a_0, a_2, a_1 ]$};

\end{tikzpicture}
\caption{Orienta{\c c}\~oes do tri\^angulo.}\label{figura12}
\end{figure}


Denotamos $[a_0, a_2, a_1] = - [a_0, a_1, a_2]$ (veja Figura \ref{figura12}). 

Uma segunda permuta{\c c}\~ao (por exemplo $[a_2, a_0, a_1]$) d\'a a volta da orienta{\c c}\~ao original.

Assim, podemos definir uma rela{\c c}\~ao de equival\^encia sobre ordens: duas ordens s\~ao equivalentes se, e somente se, podemos passar de uma \`a outra por um n\'umero par de permuta{\c c}\~oes. Podemos dar a defini{\c c}\~ao da orienta{\c c}\~ao:

{\rm A {\it orienta{\c c}\~ao} de um simplexo \'e uma classe de equival\^encia de ordens pela rela{\c c}\~ao acima.
}

\subsection{Orienta{\c c}\~ao das faces} 
Seja $\sigma = [a_0, a_1, a_2]$ um tri\^angulo orientado (pela ordem dada). 

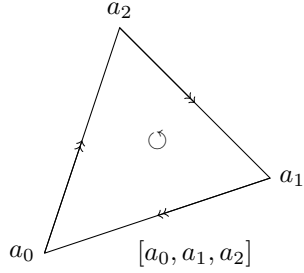
\begin{figure}[h]
\begin{tikzpicture} 
  \draw (0,0) node[left]{$a_0$}    --  (1,3) node[above]{$a_2$} -- (3, 1) node[right]{$a_1$} -- (0,0);
\node at (1.5, 1.5) {$\circlearrowleft$}; 
\node at (2,0) {$[a_0, a_1, a_2 ]$};
\draw [->>] (0,0) -- (0.5, 1.5);
\draw [->>] (1,3) -- (2,2);
\draw[->>] (3,1) -- (1.5, 0.5);
\end{tikzpicture}
\caption{Orienta{\c c}\~oes das faces do tri\^angulo.}\label{figura13}
\end{figure}


A no{\c c}\~ao intuitiva da orienta{\c c}\~ao do bordo de $\sigma$ corresponde a: 
$$[a_0, a_1] + [a_1, a_2] + [a_2, a_0], $$
(veja Figura \ref{figura13}).

Observamos que $ [a_2, a_0] = -  [a_0, a_2]$ (orienta{\c c}\~ao oposta). Assim o bordo orientado de $\sigma = [a_0, a_1, a_2]$ pode se escrever
$$ [a_1, a_2] - [ a_0, a_2] + [a_0, a_1] = [\widehat{a}_0, a_1, a_2] - [ a_0, \widehat{a}_1, a_2] + [a_0, a_1, \widehat{a}_2],$$
onde, por exemplo, $\widehat{a}_0$ significa que $a_0$ n\~ao \'e presente em $[a_1, a_2] $. Ent\~ao o bordo orientado de $\sigma = [a_0, a_1, a_2]$ \'e 
$$\sum_{i=1}^3 (-1)^i [a_0, \widehat{a}_i, a_2].$$

\subsection{Orienta{\c c}\~ao compat\'ivel}
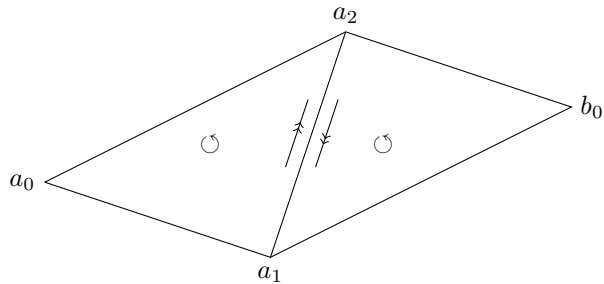
\begin{figure}[h]
\begin{tikzpicture} 
  \draw (1,1) node[left]{$a_0$}    --  (5,3) node[above]{$a_2$} -- (8, 2) node[right]{$b_0$} -- (4,0) node[below]{$a_1$}--(1,1);
\draw (4,0) -- (5,3);
 \node at (3.2,1.5) {$\circlearrowleft$};
 \node at (5.5,1.5) {$\circlearrowleft$};
\draw [->>] (4.2, 1.2) -- (4.4, 1.8);
\draw  (4.4, 1.8) --(4.5, 2.1);
\draw [->>] (4.9, 2.1) -- (4.7, 1.5);
\draw   (4.7, 1.5) -- (4.6, 1.2);

\end{tikzpicture}
\caption{Orienta{\c c}\~oes compat\'iveis.}\label{figura14}
\end{figure}

Sejam $\sigma=(a_0, a_1, a_2)$ e $\sigma' = (b_0, a_1, a_2)$ dois tri\^angulos com uma face comun $(a_1, a_2)$ (Figura \ref{figura14}).

Toda orienta{\c c}\~ao de $\sigma$, por exemplo $[a_0, a_1, a_2]$, define uma orienta{\c c}\~ao sobre $\sigma'$, que chamaremos de {\it orienta{\c c}\~ao compat\'ivel}, tal que a face comun $(a_1, a_2)$ t\^em orienta{\c c}\~oes opostas como faces de $\sigma$ e $\sigma'$. 

Esta defini{\c c}\~ao corresponde \`a no{\c c}\~ao intuitiva da compatibilidade de orienta{\c c}\~oes. 

\subsection{Orienta{\c c}\~ao de poliedros}
Seja $K$ um poliedro compacto, conexo de dimens\~ao 2. Sejam dois simplexos $\sigma$, $\tau$ de dimens\~ao 2, podemos escolher cadeias de 2-simplexos 
$$\sigma = \sigma_0, \sigma_1, \ldots, \sigma_i, \sigma_{i+1}, \ldots, \sigma_k = \tau$$
 tais que para todo $i = 0, \ldots, k-1$, os dois simplexos $\sigma_i$ e $\sigma_{i+1}$ t\^em uma face comum de dimens\~ao 1.

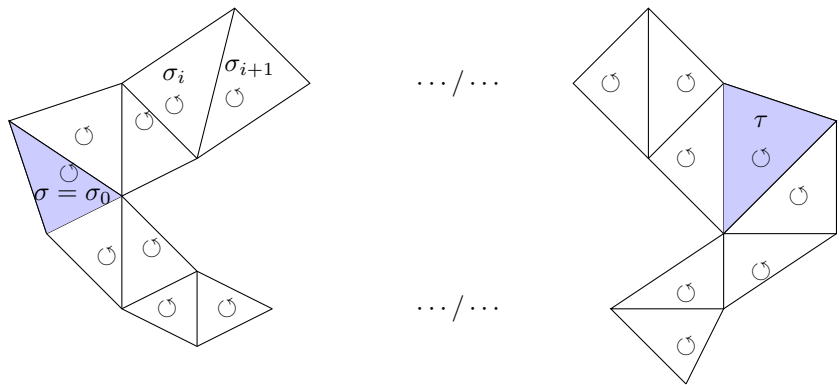
\begin{figure}[h]
\begin{tikzpicture} 
  \draw (4,1) -- (3, 0.5) -- (2,1) -- (1,2) -- (0.5, 3.5) -- (2,4) -- (3.5, 5) -- (4.5, 4) -- (3, 3) -- 
(2, 2.5) -- (3, 1.5) -- (4,1);

\draw (3, 0.5) -- (3, 1.5);
\draw (2,1) --(3, 1.5);
\draw (2,1) --(2, 2.5);
\draw (1,2) --(2, 2.5);
\draw (0.5, 3.5) --(2, 2.5);
\draw  (2,4) --(2, 2.5);
\draw  (2,4) --((3, 3);
\draw   (3.5, 5) --((3, 3);

\node at (6.5, 1) {$\cdots /  \cdots$};
\node at (6.5, 4){$\cdots /  \cdots$};

 \draw (9.5, 0) -- (8.5, 1) -- (10, 1) -- (10,2) -- (11.5,2) -- (11.5, 3.5) -- (10, 4) -- (9,5) -- (8,4) -- (9,3) -- (10, 2) -- (8.5,1);

\draw (10, 2) -- (11.5, 3.5);
\draw (10,2) -- (10,4);
\draw (9,3) -- (10,4);
\draw (9,3) -- (9,5);
\draw (9.5,0) -- (10,1);
\draw (10,1) -- (11.5, 2);

\filldraw[fill=blue!20] (1,2) -- (0.5,3.5) -- (2, 2.5);

\node at (3.4,1) {$\circlearrowleft$};
\node at (2.6,1) {$\circlearrowleft$};
\node at (2.4,1.8) {$\circlearrowleft$};
\node at (1.8,1.7) {$\circlearrowleft$};
\node at (1.3,2.8) {$\circlearrowleft$};
\node at (1.35,2.5) {$\sigma = \sigma_0$};
\node at (2.3,3.5) {$\circlearrowleft$};
\node at (1.5,3.3) {$\circlearrowleft$};

\node at (2.7,3.7) {$\circlearrowleft$};
\node at (2.7,4.1) {$\sigma_i$};
\node at (3.5,3.8) {$\circlearrowleft$};
\node at (3.7,4.2) {$\sigma_{i+1}$};

\filldraw[fill=blue!20] (10,2) -- (11.5, 3.5) -- (10,4);

\node at (9.5,0.5) {$\circlearrowleft$};
\node at (9.5,1.2) {$\circlearrowleft$};
\node at (10.5,1.5) {$\circlearrowleft$};
\node at (11,2.5) {$\circlearrowleft$};
\node at (10.5,3) {$\circlearrowleft$};
\node at (10.5,3.5) {$\tau$};
\node at (9.5,3) {$\circlearrowleft$};
\node at (9.5,4) {$\circlearrowleft$};
\node at (8.5,4) {$\circlearrowleft$};
\end{tikzpicture}
\caption{Cadeias de simplexos.}\label{figura15}
\end{figure}

\begin{figure}[h]
\hskip 20 pt
\scalebox{0.6}{\includegraphics{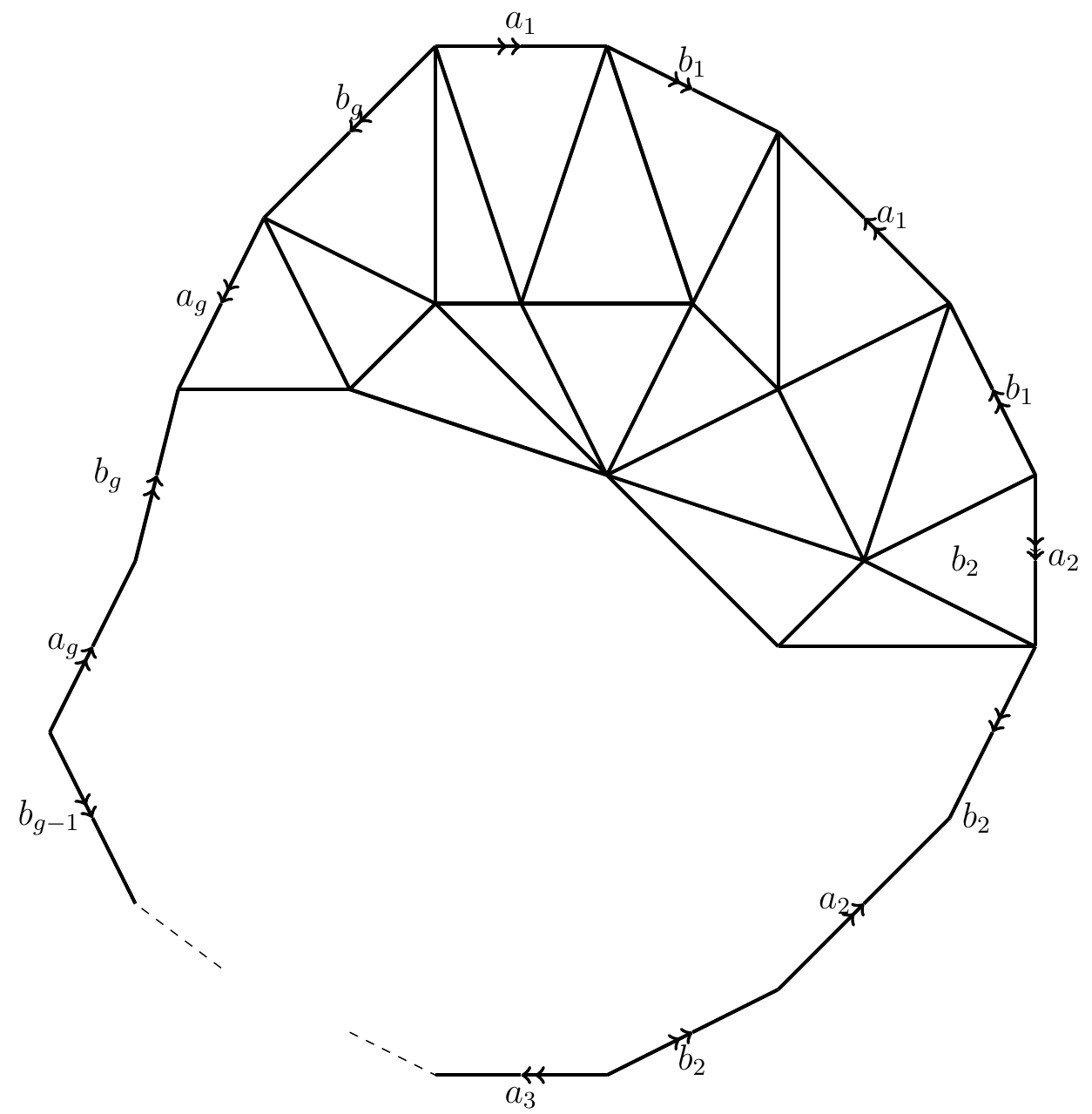}}
\caption{Representa{\c c}\~ao plana do toro de g\^enero $g$. De fato, a figura tem de ser subdivida com a observa\c c\~ao das figuras \ref{AA} e \ref{AB}, o que n\~ao fizemos aqui para ter uma figura mais leve.}\label{figura16}
\end{figure}


Escolhendo uma orienta{\c c}\~ao de $\sigma$ e uma cadeia de $\sigma$ a $\tau$, obtemos uma orienta{\c c}\~ao de $\tau$, pedindo que para todo $i = 0, \ldots, k-1$, os dois simplexos $\sigma_i$ e $\sigma_{i+1}$ t\^em orienta{\c c}\~oes compat\'iveis. 

A orienta{\c c}\~ao obtida pelo tri\^angulo $\tau$ \'e chamada de {\it orienta{\c c}\~ao induzida sobre $\tau$} a partir da orienta{\c c}\~ao de $\sigma$ ao longo da cadeia de 2-simplexos. 

\begin{definition}
Um poliedro $K$ de dimens\~ao 2, compacto, conexo \'e {\it orient\'avel} se para todo par de tri\^angulos $(\sigma, \tau)$, a orienta{\c c}\~ao induzida sobre $\tau$ a partir da orienta{\c c}\~ao de $\sigma$ n\~ao depende da cadeia de simplexos escolhida. 
\end{definition}

\subsection{Superf\'icies orient\'aveis} 

Uma superf\'icie triangulada por um poliedro $K$ \'e orient\'avel se e somente se o poliedro $K$ 
\'e orient\'avel.

As superf\'icies orient\'aveis s\~ao homeomorfas seja a esfera $\Sp$, seja ao toro $\T$, 
seja ao ``o toro de g\^enero $g$'' (Figura \ref{figura16}).

O g\^enero da esfera $\Sp$ \'e 0 : qualquer c\'irculo desenhado sobre a esfera a desconecta. 
O g\^enero do toro $\T$ \'e 1 (\'e poss\'ivel desenhar um c\'irculo no toro sem o desconectar (Figura \ref{figura15}), mas se a gente desenhar um segundo c\'irculo 
qualquer, o toro ser\'a desconectado).

\begin{figure}[h]
\hskip 20 pt
\scalebox{0.60}{\includegraphics{AD.jpg}}
\caption{O toro \'e  de g\^enero 1 }\label{figura15}
\end{figure}

A caracter\'istica de Euler-Poincar\'e de uma superf\'icie orient\'avel de g\^enero $g$ vale $2 -2g$ (Simon Antoine-Jean  Lhuillier \cite{Lhu}, 1812). 

O toro de g\^enero $g$ \'e a soma conexa de um toro de g\^enero $g-1$ e do toro. 

Uma superf\'icie compacta sem bodo \'e orient\'avel se e somente se ela pode ser mergulhada em $\R^3$.

\subsection{Superf\'icies n\~ao orient\'aveis}
As superf\'icies n\~ao orient\'aveis bem conhecidas s\~ao a {\it garrafa de Klein} e o espa{\c c}o projeto ${\mathbb{P}}^2$. 

O espa{\c c}o projetivo ${\mathbb{P}}^2$ \'e o conjunto de todas retas do espa{\c c}o euclidiano $\R^3$ passando pela origem. 
Uma maneira f\'acil para representar o espa{\c c}o projetivo \'e de considerar em $\R^3$ a esfera $\Sp^2$ centrada na origem e de raio 1.

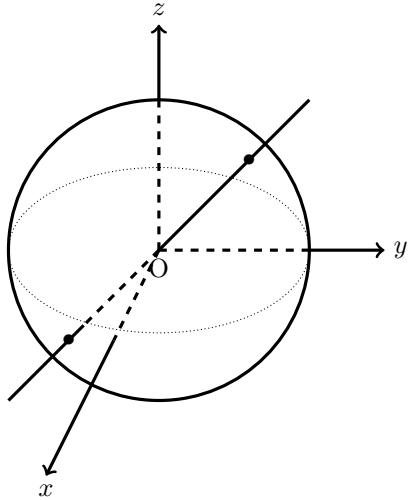
\begin{figure}[h]
\begin{tikzpicture} 
\draw[very thick] (2,2) circle (2cm);
\draw[densely dotted] (2,2) ellipse (2cm and 1.1cm);
\draw[very thick, dashed] (2,2) --(2,4);
\draw[very thick, ->] (2,4) -- (2,5);
\node at (2, 5)[above] {$z$};

\draw[very thick, dashed] (2,2) -- (4,2);
\draw[very thick, ->] (4,2) -- (5,2);
\node at (5,2)[right] {$y$};

\draw[very thick] (2,2) -- (4,4);
\node at (3.2,3.2) {$\bullet$};
\draw[very thick, dashed] (2,2) -- (1,1);
\draw[very thick] (1,1) -- (0,0);
\node at (0.8,0.8) {$\bullet$};

\draw[very thick, dashed] (2,2) -- (1.4, 0.8);
\draw[->, very thick] (1.4, 0.8) -- (0.5, -1);

\node at  (0.5, -1)[below] {$x$};

\node at (2,2)[below] {O};

\end{tikzpicture} 
\caption{Representa{\c c}\~ao do espa\c co projetivo ${\mathbb{P}}^2$.}\label{figura17}
\end{figure}[h]


Cada reta de $\R^3$ passando pela origem  (e fora do plano $0xy$) encontra a semiesfera norte em um ponto $x_0$.
 Este ponto \'e um representante da reta. 

\begin{figure}[h]
\hskip 20 pt
\scalebox{0.7}{\includegraphics{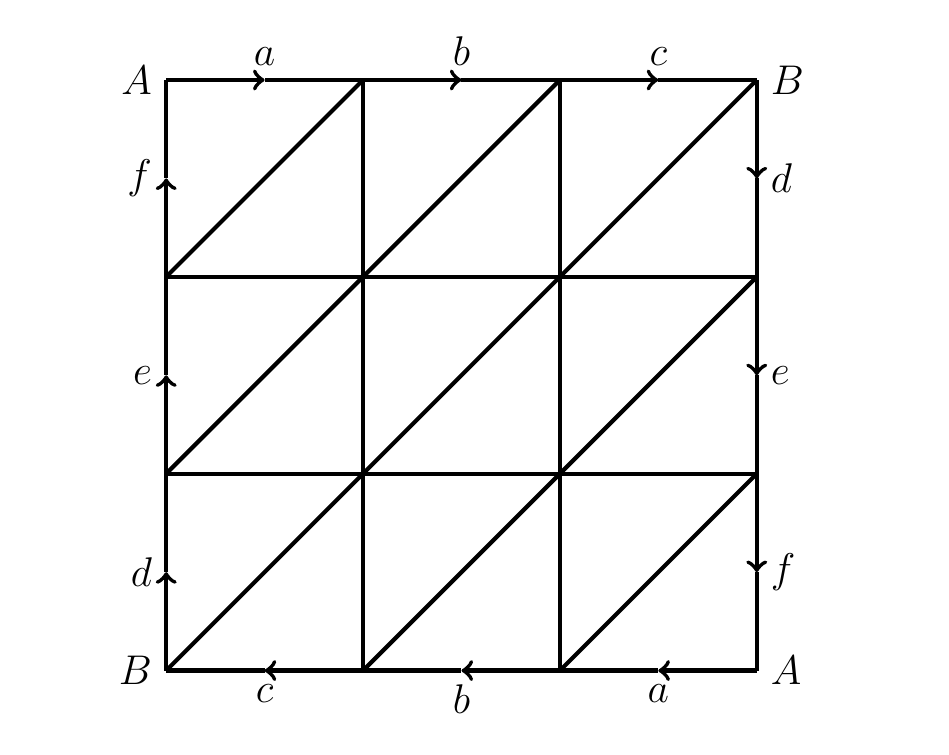}}
\caption{Representa{\c c}\~ao plana do espa\c co projetivo ${\mathbb{P}}^2$.}\label{Projectif}
\end{figure}

As retas do plano $0xy$ encontram a semiesfera norte em dois pontos diametralmente opostos, que temos que identificar. Assim obtemos uma representa{\c c}\~ao do espa{\c c}o projetivo ${\mathbb{P}}^2$, 
como sendo a semiesfera norte com identifica\c c\~ao de pontos diametralmente opostos no seu bordo. 
Outra maneira \'e de considerar a representa{\c c}\~ao plana da figura \ref{Projectif} 
que \'e, obviamente uma representa{\c c}\~ao \'equivalente de ${\mathbb{P}}^2$.

A caracter\'istica de Euler-Poincar\'e de ${\mathbb{P}}^2$ vale 
$$\chi({\mathbb{P}}^2) = n_0 - n_1 + n_2 = 10 - 27 + 18 = +1.$$
O gen\^ero de ${\mathbb{P}}^2$ vale 1. 

\begin{figure}[h]
\hskip 20 pt
\scalebox{0.70}{\includegraphics{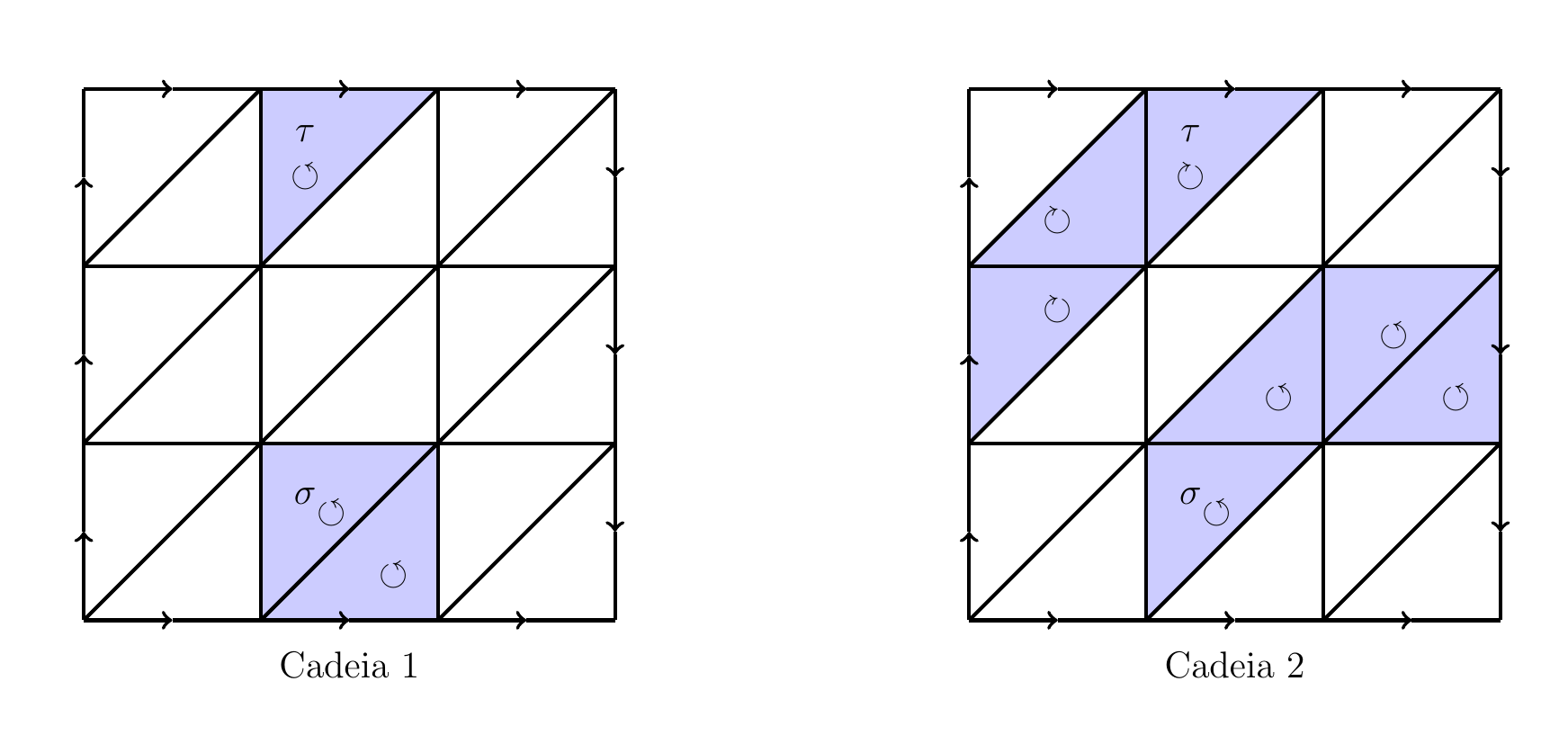}}
\caption{O espa\c co projetivo ${\mathbb{P}}^2$ n\~ao \'e orient\'avel: 
Indo de $\sigma$ at\'e $\tau$, pela cadeia de simplexos 1 e pela cadeia de simplexos 2, 
d\'a orienta{\c c}\~oes diferentes de $\tau$.}\label{figura18}
\end{figure}

A figura \ref{figura18} illustra o fato que ${\mathbb{P}}^2$ n\~ao \'e orient\'avel. 

\begin{figure}[h]
\hskip 20 pt
\scalebox{0.70}{\includegraphics{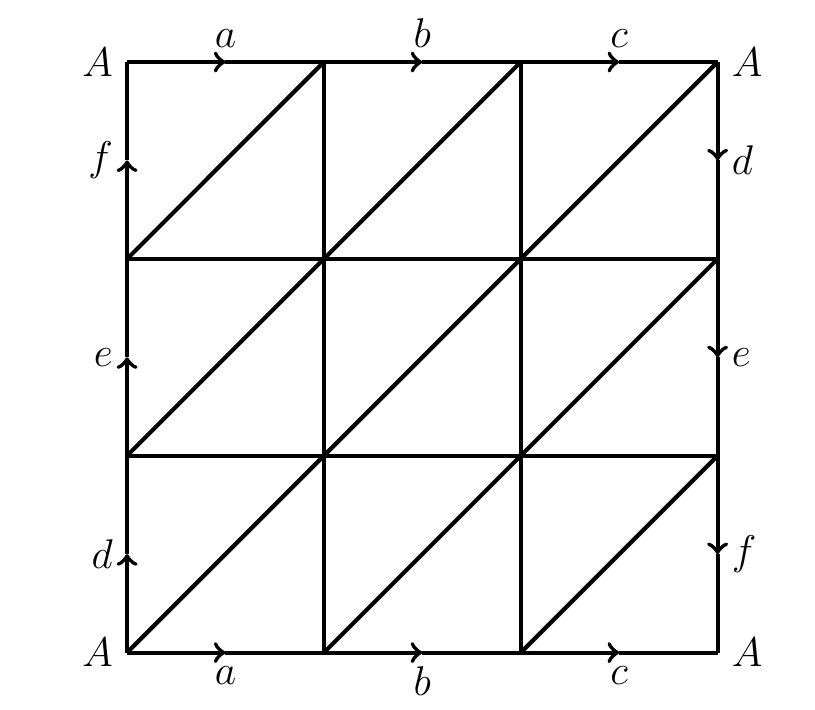}}
\caption{Representa{\c c}\~ao plana da garrafa de Klein.}\label{Klein}
\end{figure}

A garrafa de Klein $G$ \'e a superf\'icie representada de uma maneira plana como na figura \ref{Klein}.

A garrafa de Klein n\~ao \'e orient\'avel, a sua caracter\'istica de Euler-Poincar\'e \'e igual a
$$\chi(G) = n_0 - n_1 + n_2 = 9 -27 + 18 = 0.$$

\begin{figure}[h]
\hskip 20 pt
\scalebox{0.40}{\includegraphics{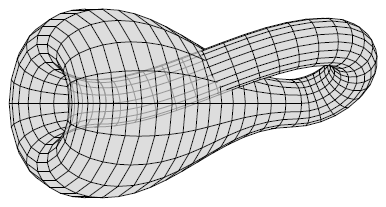}}
\caption{A garrafa de Klein.}\label{Klein0}
\end{figure}

No espa\c co euclidiano, n\~ao \'e poss\'ivel representar a garrafa de Klein sem ``autointersec\c c\~ao"
(veja Figura \ref{Klein0}). E imposs\'ivel preencher de \'agua (ou outro l\'iquido !)  na garrafa de Klein.

A garrafa de Klein \'e a soma conexa de dois planos projetivos.

De uma maneira geral, para uma superf\'icie n\~ao orient\'avel, de g\^enero $k$, a caracter\'istica de Euler-Poincar\'e vale $2 -k$. Uma superf\'icie n\~ao orient\'avel de g\^enero $k$ \'e homeomorfa \`a soma conexa de $k$ espa{\c c}os  ${\mathbb{P}}^2$.

\subsection{Generaliza{\c c}\~oes}
\subsubsection{Superf\'icies singulares}
\begin{figure}[h]
\hskip 20 pt
\scalebox{1.40}{\includegraphics{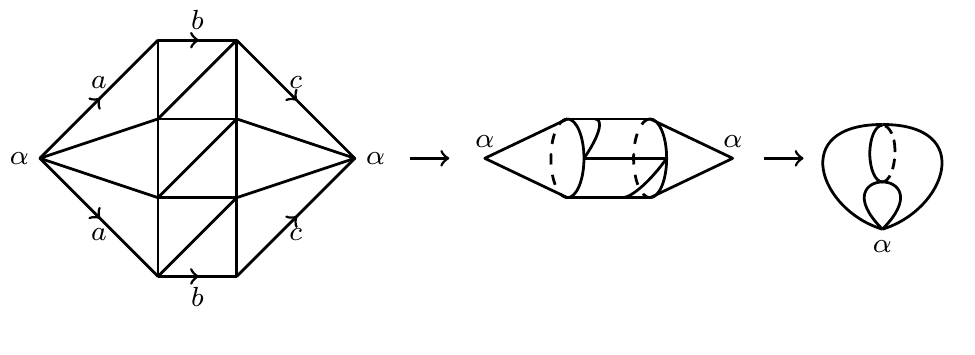}}
\caption{O toro pin{\c c}ado.}\label{toropincado}
\end{figure}

\'E poss\'ivel de calcular a caracter\'istica de Euler-Poincar\'e de uma superf\'icie singular do mesmo jeito do que no caso liso. Por exemplo, a caracter\'istica de Euler-Poincar\'e do toro pin{\c c}ado vale $1$ 
qualquer que seja a tri\^angula{\c c}\~ao do toro pin{\c c}ado. Na figura \ref{toropincado} temos
$$\chi(X) = n_0 - n_1 + n_2 = 7 -18 + 12 = 1.$$

\subsubsection{Dimens\~oes maiores} 
No caso de espa{\c c}os de dimens\~ao maior do que 2, 
consideramos espa{\c c}os $X$ compactos, conexos, de dimens\~ao pura $k$. 
Isto \'e todo ponto de $X$ pertence ao fecho de (ao menos) um simplexo de dimens\~ao $k$ de uma tri\^angula{\c c}\~ao $K$ de $X$. 

\'E claro que, neste caso, $K$ admite simplexos de todas dimens\~oes de $0$ at\'e $k$. 
Um simplexo padr\~ao de dimens\~ao $i$ \'e representado no espa\c co euclidiano 
$\R^{i+1}$ de coordenadas $(x_1, \ldots, x_{i+1})$ como o conjunto de pontos tais que $\sum_{\alpha=1}^{i+1} x_\alpha =1$ e $0\le i \le 1$ para tudo $i$.

Denotamos por $n_i$ o n\'umero de simplexos de dimens\~ao $i$ em $K$ ($i = 0, \ldots, k$).  

Poincar\'e mostrou (no ano 1893) que a soma 
$$\chi(X) = \sum_{i=0}^k (-1)^i n_i$$
n\~ao depende da tri\^angula{\c c}\~ao $K$ do espa{\c c}o $X$. Este \'e a defini{\c c}\~ao da  caracter\'istica de Euler-Poincar\'e de um espa{\c c}o $X$ de dimens\~ao $k$, que definitivamente deveria ser chamada de 
{\it caracter\'istica de Descartes-Poincar\'e.}

\end{document}